\documentclass[12pt]{amsart}
\makeatletter \@addtoreset{equation}{section}
\makeatother

\usepackage{amssymb}
\usepackage{mathrsfs}
\newcommand{\Pic}{\operatorname{Pic}}
\newcommand{\Bs}{\operatorname{Bs}}
\newcommand{\Cr}{\operatorname{Cr}}
\newcommand{\Aut}{\operatorname{Aut}}
\newcommand{\Bir}{\operatorname{Bir}}
\newcommand{\Sing}{\operatorname{Sing}}

\newcommand{\wt}{\operatorname{wt}}

\newcommand{\PP}{\mathbb{P}}
\newcommand{\ZZ}{\mathbb{Z}}
\newcommand{\CC}{\mathbb{C}}
\newcommand{\QQ}{\mathbb{Q}}
\newcommand{\EEE}{{\mathscr{E}}}

\newcommand{\PGL}{{\operatorname{PGL}}}

\newcommand{\GL}{{\operatorname{GL}}}
\newcommand{\rk}{\operatorname{rk}}

\newcommand{\xref}[1]{{\rm \ref{#1}}}

\newcommand{\muu}{{\boldsymbol{\mu}}}

\newtheorem{theorem}[equation]{Theorem}
\newtheorem{proposition}[equation]{Proposition}
\newtheorem{proposition-definition}[equation]{Proposition-Definition}
\newtheorem{lemma}[equation]{Lemma}
\newtheorem{corollary}[equation]{Corollary}

\newtheorem{corollary-notation}[equation]{Corollary-Notation}

\theoremstyle{definition}

\newtheorem{de}[equation]{}
\newtheorem{definition}[equation]{Definition}

\newtheorem{pusto}[equation]{}
\newtheorem{remark}[equation]{Remark}

\theoremstyle{remark}

\date{}
\title[]{$p$-elementary
subgroups of the Cremona group of rank $3$}
\author[]{Yuri Prokhorov}
\address{Yuri Prokhorov: Department 
of Algebra, Faculty of Mathematics, Moscow State
 University, Moscow 117234, Russia }
 \email{prokhoro@gmail.com}
\thanks{The author was partially supported by RFBR, grant
No. 08-01-00395-a and Leading Scientific
Schools (grants No. 1983.2008.1, 1987.2008.1).
}

\begin{document}
\begin{abstract}
For the subgroups of the Cremona group $\mathrm{Cr}_3(\mathbb C)$
having the form $(\boldsymbol{\mu}_p)^s$, where $p$ is prime, 
we obtain an upper bound for $s$. 
Our bound is sharp if $p\ge 17$. 
\end{abstract}
\maketitle
\section{Introduction}
Let $\Bbbk$ be an algebraically closed field. The {\it Cremona group}
$\Cr_{n}(\Bbbk)$ is the group of birational transformations of $\PP^n_{\Bbbk}$, or
equivalently the group of $\Bbbk$-automorphisms of the field
$\Bbbk(x_1,\dots,x_n)$. 
Finite subgroups of $\Cr_{2}(\CC)$ are completely classified (see \cite{Dolgachev-Iskovskikh}
and references therein).
In contrast, subgroups of $\Cr_{n}(\Bbbk)$ for $n\ge 3$ 
are not studied well (cf. \cite{Prokhorov2009b}). 

In the present paper we study some kind of abelian subgroups of $\Cr_3(\CC)$.
Let $p$ be a prime number.
We say that a group $G$ is \textit{$p$-elementary} if 
$G\simeq (\muu_p)^s$ for some positive integer $s$.
In this case $s$ is called the \textit{rank} of $G$ and denoted by 
$\rk G$.

\begin{theorem}[{\cite{Beauville2007}}]
\label{theorem-cr2}
Let $p$ be a prime $\neq \operatorname{char}(\Bbbk)$
and let $G\subset \Cr_2(\Bbbk)$ be a $p$-elementary subgroup.
Then:
\[
\rk G \le 2+\delta_{p,3}+2\delta_{p,2}
\]
where $\delta_{i,j}$ is Kronecker's delta.
Moreover, for any such $p$ 
this bound is attained for some $G$. These ``maximal'' groups $G$ are classified up
to conjugacy in $\Cr_2(\Bbbk)$.
\end{theorem}

More generally, instead of $\Cr_n(\Bbbk)$  we also can consider 
the group $\Bir(X)$ of birational automorphisms 
of an arbitrary rationally connected variety $X$.
Our main result is the following

\begin{theorem}
\label{theorem-cr3}
Let $X$ be a rationally connected threefold defined 
over a field of characteristic $0$, let $p$ be a prime,
and let $G\subset \Bir (X)$ be a $p$-elementary subgroup.
Then
\begin{equation}
\label{eq-main}
\rk G\le
\begin{cases}
7&\text{if $p=2$,}
\\
5&\text{if $p=3$,}
\\
4&\text{if $p=5,\, 7,\,11$, or $13$,}
\\
3&\text{if $p\ge 17$.}
\end{cases}
\end{equation}
For any prime $p\ge 17$ this bound is attained for some subgroup $G\subset \Cr_3(\CC)$.
\textup(However we do not assert that the bound \eqref{eq-main}
is sharp for $p\le 13$\textup).
\end{theorem}

\begin{remark}
(i)
Note that $\Cr_1(\Bbbk)\simeq \PGL_2(\Bbbk)$.
Hence for any prime $p$ $\neq \operatorname{char}(\Bbbk)$
and any $p$-elementary subgroup $G\subset \Cr_1(\Bbbk)$, 
we have $\rk G\le 1+ \delta_{p,2}$
(see, e.g., \cite[Lemma 2.1]{Beauville2007}).

(ii)
Since $\Cr_1(\Bbbk)\times \Cr_2(\Bbbk)$ admits (a lot of) embeddings to 
$\Cr_3(\Bbbk)$, the group $\Cr_3(\Bbbk)$ contains a $p$-elementary 
subgroup $G$ of rank $3+\delta_{p,3}+3\delta_{p,2}$.
This shows the last assertion of our theorem.
\end{remark}

The following consequence of 
Theorem \ref{theorem-cr3}
was proposed by A. Beauville.
\begin{corollary}
\label{cor-cr}
The group $\Cr_3(\CC)$ is not isomorphic to $\Cr_n(\CC)$ for $n\neq 3$ as an abstract group. 
\end{corollary}

\begin{proof}
Denote by $\xi(n,p)$ the maximal rank of a $p$-elementary group contained in $\Cr_n(\CC)$.
Then $\xi(2,17)=2<\xi(3,17)=3$ and $\xi(n,17)\ge n$ by Theorems \ref{theorem-cr2} and \ref{theorem-cr3}.
\end{proof}

Our method is a generalization of the method
used for study of finite subgroups of $\Cr_2(\Bbbk)$ 
\cite{Beauville2007}, \cite{Dolgachev-Iskovskikh}. Similarly to \cite{Prokhorov2009b}
we use the equivariant three-dimensional minimal model program. 
This way easily allows us to reduce the problem to the
study of automorphism groups of some (not necessarily smooth) Fano threefolds.

\textbf{Acknowledgments.}
I would like  to thank J.-P. Serre for asking me the question considered here
and  for useful comments.
I am also grateful to A. Beauville for 
proposing me Corollary \ref{cor-cr} and for his interest on my paper.
This paper was written
at the IHES (Bures-sur-Yvette) during my visit in 2009. 
I am grateful to IHES for the support and hospitality.
Finally I would like to thank the referee for several suggestions that make
paper more readable and clear.

\section{Preliminaries}
Clearly, we may assume that $\Bbbk=\CC$. 
All the groups in this paper are multiplicative.
In particular, we denote a cyclic group of order $n$ by $\muu_n$.


\begin{de}\label{Terminal-singularities}
\textbf{Terminal singularities.}
We need a few facts on the classification of 
three-dimensional 
terminal singularities (see  \cite{Mori-1985-cla}, \cite{Reid-YPG1987}).
Let $(X\ni P)$ be a germ of a three-dimensional terminal singularity.
Then $(X\ni P)$ is isolated, i.e., $\Sing(X)=\{P\}$.
The \textit{index} of $(X\ni P)$ is the minimal positive integer
$r$ such that $rK_X$ is Cartier.
If $r=1$, then $(X\ni P)$ is Gorenstein. In this case $(X\ni P)$ is analytically 
isomorphic to a hypersurface 
singularity in $\CC^4$ of multiplicity $2$.
Moreover, any Weil $\QQ$-Cartier divisor $D$ on $(X\ni P)$ is Cartier.
If $r>1$, then there is a cyclic, \'etale outside of $P$ 
cover $\pi : (X^{\sharp}\ni P^{\sharp})\to (X\ni P)$ 
of degree $r$ such that $(X^{\sharp}\ni P^{\sharp})$ is a Gorenstein 
terminal singularity (or a smooth point). This $\pi$ is called
the \textit{index-one cover} of $(X\ni P)$. 
\end{de}
\begin{theorem}[\cite{Mori-1985-cla}, \cite{Reid-YPG1987}]
\label{terminal-lemma}
In the above notation $(X^{\sharp}\ni P^{\sharp})$ is 
analytically $\muu_r$-isomorphic to a hypersurface 
in $\CC^4$ with $\muu_r$-semi-invariant
\footnote{In invariant theory people often say ``relative invariant'' rather then 
``semi-invariant''. We prefer to use the terminology of \cite{Mori-1985-cla}.} coordinates
$x_1,\dots,x_4$,
and the action is given by 
\[
 (x_1,\dots,x_4) \longmapsto (\varepsilon^{a_1} x_1,\dots,\varepsilon^{a_4} x_4)
\]
for some primitive $r$-th root of unity $\varepsilon$, where
one of the following holds:
\begin{enumerate}
 \item 
$(a_1,\dots,a_4)\equiv (1,-1,a_2,0)\mod r$, \quad
$\gcd(a_2,r)=1$, 
 \item 
$r=4$ and $(a_1,\dots,a_4)\equiv (1,-1,1,2)\mod 4$.
\end{enumerate}
\end{theorem}

\begin{definition}
A \emph{$G$-variety} is a variety $X$ provided with 
a biregular faithful action of a finite group $G$.
We say that a normal $G$-variety 
$X$ is \emph{$G\QQ$-factorial} if any $G$-invariant Weil 
divisor on $X$ is $\QQ$-Cartier. 
A projective normal $G$-variety $X$ is called \emph{$G\QQ$-Fano}
if it is $G\QQ$-factorial, has at worst terminal 
singularities, $-K_X$ is ample, and $\rk \Pic (X)^G=1$. 
\end{definition}

\begin{lemma}
\label{lemma-fixed-point}
Let $(X\ni P)$ be a germ of a threefold terminal singularity 
and let $G\subset \Aut(X\ni P)$ be a $p$-elementary subgroup. 
Then $\rk G\le 3+\delta_{2,p}$.
\end{lemma}
\begin{proof}
Assume that $\rk G\ge   4+\delta_{2,p}$.
First we consider the case where $(X\ni P)$ is Gorenstein.
The group $G$ acts faithfully on
the Zariski tangent space $T_{P,X}$, so 
$G\subset \GL (T_{P,X})$, where $\dim T_{P,X}=3$ or $4$.
If $\dim T_{P,X}=3$, then 
$G$ is contained in a maximal torus of $\GL_3(\CC)$, so 
$\rk G\le 3$ and we are done. 
Thus we may assume that $\dim T_{P,X}=4$. 
Take semi-invariant coordinates 
$x_1,\dots, x_4$ in $T_{P,X}$. There is a $G$-equivariant 
analytic embedding $(X\ni P)\subset \CC_{x_1,\dots, x_4}^4$.
As above, $\rk G\le 4$. Thus we may assume that $\rk G\le 4$ and $p>2$.
Let $\phi(x_1,\dots, x_4)=0$ be an equation of $X$, where $\phi$ is a $G$-semi-invariant function. 
Regard $\phi$ as a power series and write $\phi=\sum_d \phi_d$,
where $\phi_d$ is the sum of all monomials of degree $d$.
Since the action of $G$ on $x_1,\dots, x_4$ is linear, all the $\phi_d$'s 
are semi-invariants of the same $G$-weight $w=\wt \phi_d$. 
Hence, for any $\phi_d$, $\phi_{d'}\neq 0$ we have $d-d'\equiv 0 \mod p$.
Since $(X\ni P)$ is a terminal singularity, $\phi_2\neq 0$ and so $\phi_3=0$.
Recall that $G\simeq (\muu_p)^4$, $p\ge 3$.
In this case, $\phi_2$ must be a monomial.
Thus up to permutations of coordinates and scalar multiplication
we get either $\phi_2=x_1^2$ or $\phi_2=x_1x_2$.
In particular, we have $\rk \phi_2\le 2$ and $\phi_3=0$.
This contradicts 
the classification of terminal singularities \cite{Mori-1985-cla}, \cite{Reid-YPG1987}.

Now assume that $(X\ni P)$ is non-Gorenstein of index $r>1$.
Consider the index-one cover $\pi\colon (X^\sharp\ni P^\sharp)\to 
(X\ni P)$ (see \ref{Terminal-singularities}). 
Here $(X^\sharp\ni P^\sharp)$ is a Gorenstein terminal 
point and the map $X^\sharp\setminus \{P^\sharp\}\to 
X\setminus \{P\}$ can be regarded as the topological universal cover. 
Hence there exists a natural lifting
$G^\sharp\subset \Aut (X^\sharp\ni P^\sharp)$ fitting to the following exact sequence
\begin{equation}
\label{eq-seq-1}
1\longrightarrow \muu_r \longrightarrow G^\sharp \longrightarrow G \longrightarrow 1.
\end{equation}
It is sufficient to show that there exists a 
subgroup $G^{\bullet}\subset G^\sharp$ isomorphic to $G$
(but we do not assert that the sequence splits).
Indeed, in this case $G^{\bullet}\simeq G$ acts faithfully on the terminal Gorenstein
singularity $(X^\sharp\ni P^\sharp)$ and we can apply the above 
considered case. 
We may assume that $G^\sharp$ is not abelian 
(otherwise a subgroup $G^{\bullet}\simeq G$ obviously exists). 
The group $G^\sharp$ permutes 
eigenspaces of $\muu_{r}$. 
By Theorem \ref{terminal-lemma} the subspace $T:=\{x_4=0\}\subset \CC^4_{x_1,\dots,x_4}$
is $G^\sharp$-invariant and $\muu_{r}$ acts on any eigenspace $T_1\subset T$ faithfully.
On the other hand, by \eqref{eq-seq-1} we see that the derived subgroup
$[G^\sharp, G^\sharp]$ is contained in $\muu_{r}$.
In particular, $[G^\sharp, G^\sharp]$ 
is abelian and also acts on any eigenspace $T_1\subset  T$ faithfully.
Since $\dim T=3$, this implies that the representation of 
$G^\sharp$ on $T$ is irreducible
(otherwise $T$ has a one-dimensional subrepresentation, say $T_1$,
and the kernel of the map $G\to \GL(T_1)\simeq \CC^*$
must contain $[G^\sharp, G^\sharp]$). 
Hence eigenspaces of $\muu_{r}$ have the same dimension and so
$\muu_{r}$ acts on $T$ by scalar multiplication.
By Theorem \ref{terminal-lemma} this is possible only if $r=2$.

Let $G^\sharp_p\subset G^\sharp$ be a Sylow $p$-subgroup.
If $\muu_r\cap G^\sharp_p=\{1\}$, then $G^\sharp_p\simeq G$ and we are done.
Thus we assume that $\muu_r\subset G^\sharp_p$, so $p=r=2$ and $G^\sharp_p=G^\sharp$.
But then $G^\sharp$ is a $2$-group, so the dimension of its irreducible 
representation must be a power of $2$. Hence  
$T$ is reducible, a contradiction.
\end{proof}

\begin{lemma}
\label{lemma-fixed-points}
Let $X$ be a $G$-threefold with isolated singularities.
\begin{enumerate}
 \item 
If there is a curve $C\subset X$ of $G$-fixed points,
then $\rk G\le 2$. 
 \item 
If there is a surface $S\subset X$ of $G$-fixed points,
then $\rk G\le 1$. If moreover $S$ is singular along a curve,
then $G=\{1\}$.
\end{enumerate}
\end{lemma}
\begin{proof}[Sketch of the Proof]
Consider the action of $G$ on 
the tangent space to $X$ at a general point of $C$ (resp. $S$). 
\end{proof}

\textbf{$G$-equivariant minimal model program.}
Let $X$ be a rationally connected three-dimensional algebraic variety
and let $G\subset \Bir(X)$ be a finite subgroup.
By shrinking $X$ we may assume that $G$ acts on $X$ biregularly.
The quotient $Y=X/G$ is quasiprojective, so 
there exists a projective completion $\hat Y\supset Y$.
Let $\hat X$ be the normalization of $\hat Y$ in the function field
$\CC(X)$. Then $\hat X$ is a projective variety birational to $X$ 
admitting a biregular action of $G$.
There is an equivariant resolution of singularities $\tilde X\to \hat X$,
see \cite{Abramovich-Wang}.
Run the $G$-equivariant minimal model program: $\tilde X\to \bar X$,
see \cite[0.3.14]{Mori-1988}.
Running this program we stay in the category of projective normal varieties 
with at worst terminal $G\QQ$-factorial singularities. 
Since $X$ is rationally connected, on the final step we get 
a Fano-Mori fibration $f: \bar X\to Z$. Here $\dim Z<\dim X$, $Z$ is normal,
$f$ has connected fibers, the anticanonical Weil
divisor $-K_{\bar X}$ is ample over $Z$, and the relative $G$-invariant Picard 
number $\rho(\bar X)^G$ is one.
Obviously, we have the following possibilities:
\begin{enumerate}
\item 
$Z$ is a rational surface and a general fiber $F=f^{-1}(y)$ is a conic;
\item
$Z\simeq \PP^1$ and a general fiber $F=f^{-1}(y)$ is a smooth del Pezzo surface;
\item
$Z$ is a point and $\bar X$ is a $G\QQ$-Fano threefold.
\end{enumerate}

\begin{proposition}
In the above notation assume that $Z$ is not a point.
Then $\rk G\le 3+\delta_{p,3}+3\delta_{p,2}$.
In particular, \eqref{eq-main} holds.
\end{proposition}

\begin{proof}
Let $G_0\subset G$ be the kernel of the homomorphism 
$G\to \Aut(Z)$. The group $G_1:=G/G_0$ acts effectively on $Z$ 
and $G_0$ acts effectively on a 
general fiber $F\subset X$ of $f$. Hence, $G_1\subset \Aut(Z)$ and $G_0\subset \Aut(F)$.
Clearly, $G_0$ and $G_1$ are $p$-elementary groups with 
$\rk G_0 +\rk G_1=\rk G$.
Assume that $Z\simeq \PP^1$. Then $\rk G_1\le 1+\delta_{p,2}$.
By Theorem \ref{theorem-cr2} we obtain
$\rk G_0\le 2+\delta_{p,3}+2\delta_{p,2}$.
This proves 
our assertion in the case $Z\simeq \PP^1$.
The case $\dim Z=2$ is treated similarly.
\end{proof}

\begin{de}
\label{Main-assumption}
\textbf{Main assumption.}
Thus from now on we assume that we are in the case (iii).
Replacing $X$ with $\bar X$ we may assume that our original 
$X$ is a $G\QQ$-Fano threefold.

The group $G$ acts naturally on $H^0(X,-K_X)$.
If $H^0(X,-K_X)\neq 0$, then there exists 
a $G$-semi-invariant section $s\in H^0(X,-K_X)$
(because $G$ is an abelian group).
This section gives us an invariant member
$S\in |{-}K_X|$.
\end{de}

\begin{lemma}
\label{lemma-hyp-sect-0}
Let $X$ be a $G\QQ$-Fano threefold, where $G$ is a $p$-elementary
group with $\rk G\ge \delta_{p,2}+4$. 
Let $S$ be an invariant Weil divisor 
such that $-(K_X+S)$ is nef.  
Then the pair $(X,S)$ is log canonical \textup(LC\textup).
\end{lemma}
\begin{proof}
Assume that the pair $(X,S)$ is not LC. 
Since $S$ is $G$-invariant and $\rho(X)^G=1$, we see that $S$ is 
numerically proportional to $K_X$. Since  $-(K_X+S)$ is nef, $S$ is ample.
We apply quite standard
connectedness arguments of Shokurov \cite{Shokurov-1992-e-ba}
(see, e.g., \cite[Prop. 2.6]{Mori-Prokhorov-2008d}): 
for a suitable $G$-invariant boundary $D$, the pair 
$(X,D)$ is LC, the divisor $-(K_X+D)$ is ample,
and the minimal locus $V$ of log canonical singularities is also $G$-invariant. 
Moreover, $V$ is either a point or a smooth rational curve.
By Lemma \ref{lemma-fixed-point} we may assume that $G$ has no fixed points.
Hence, $V\simeq \PP^1$ and we have a map 
$\varsigma: G\to \Aut(\PP^1)$. If $p>2$, then $\varsigma(G)$ is a cyclic group, 
so $G$ has a fixed point, a contradiction.
Let $p=2$ and let $G_0=\ker \varsigma$. 
By Lemma \ref{lemma-fixed-points} $\rk G_0\le 2$.
Therefore $\rk \varsigma(G_0)\ge 3$.
Again we get a contradiction. 
\end{proof}

\begin{lemma}
\label{lemma-hyp-sect}
Let $X$ be a $G\QQ$-Fano threefold, where $G$ is a $p$-elementary
group with 
\begin{equation}
\label{main-assumption-3}
\rk G\ge 
\begin{cases}
7&\quad \text{if $p=2$,}
\\
5&\quad \text{if $p=3$,}
\\
4&\quad \text{if $p\ge 5$.}
\\
\end{cases}
\end{equation}
Let $S\in |{-}K_X|$ be a $G$-invariant member. 
Then we have
\begin{enumerate}
\item 
Any component $S_i\subset S$ is either rational or 
birationally ruled over an elliptic curve.
\item
The group $G$  acts transitively on the components of $S$.
\item
For the stabilizer $G_{S_i}$ we have $\rk G_{S_i}\le \delta_{p,2}+4$. 
\item
The surface $S$ is reducible \textup(and reduced\textup).
\end{enumerate}
\end{lemma}

\begin{proof}
By Lemma \ref{lemma-hyp-sect-0} the pair $(X,S)$ is LC.
Assume that $S$ is normal (and irreducible).
By the adjunction formula $K_S\sim 0$.
We claim that $S$ has at worst Du Val singularities. Indeed, otherwise 
by the Connectedness Principle \cite[Th. 6.9]{Shokurov-1992-e-ba}
$S$ has at most two non-Du Val points. 
If $p>2$, these points must be $G$-fixed. This  contradicts Lemma \ref{lemma-fixed-point}.
Otherwise $p=2$ and these points are fixed for an index two subgroup 
$G^{\bullet}\subset G$. Again we get a contradiction by Lemma \ref{lemma-fixed-point}.
Thus we may assume that $S$ has at worst Du Val singularities.
Let $\Gamma$ be the image of $G$ in $\Aut(S)$.
By Lemma \ref{lemma-fixed-points} $\rk G\le \rk \Gamma+1$.
Let $\tilde S\to S$ be the minimal resolution.
Here $\tilde S$ is a smooth K3 surface.
The natural representation of $\Gamma$ on $H^{2,0}(\tilde S)$ 
induces the exact sequence (see \cite{Nikulin-1980-aut})
\[
1 \longrightarrow \Gamma_0 \longrightarrow \Gamma \longrightarrow \Gamma_1 \longrightarrow 1,
\]
where $\Gamma_0$ (resp. $\Gamma_1$) is the kernel (resp. image) 
of the representation of $\Gamma$ on $H^{2,0}(\tilde S)$.
The group $\Gamma_1$ is cyclic. Hence either $\Gamma_1=\{1\}$ or 
$\Gamma_1\simeq \muu_p$. In the second case by \cite[Cor. 3.2]{Nikulin-1980-aut}
$p\le 19$. Further, according to \cite[Th. 4.5]{Nikulin-1980-aut} we have 
\[
\rk \Gamma_0\le 
\begin{cases}
4 & \text{if $p=2$}
\\
2 & \text{if $p=3$}
\\
1 & \text{if $p=5$ or $7$}
\\
0& \text{if $p>7$.}
\end{cases} 
\]
Combining this we obtain a contradiction with 
\eqref{main-assumption-3}.

Now assume that $S$ is not normal.
Let $S_i\subset S$ be an irreducible component (the case $S_i=S$
is not excluded).
Let $\nu\colon S'\to S_i$
be the normalization. 
Write $0\sim \nu^*(K_X+S_i)=K_{S'}+D'$,
where $D'$ is the different, see \cite[\S 3]{Shokurov-1992-e-ba}. 
Here $D'$ is an effective reduced divisor
and the pair is LC \cite[3.2]{Shokurov-1992-e-ba}.
Since $S$ is not normal, $D'\neq 0$.
Consider the minimal resolution 
$\mu\colon \tilde S\to S'$ and let $\tilde D$ be the crepant 
pull-back of $D'$, that is,
$\mu_*\tilde D=D'$ and
\[
K_{\tilde S}+\tilde D=\mu^*(K_{S'}+D')\sim 0.
\]
Here $\tilde D$ is again an effective reduced divisor.
Hence $\tilde S$ is a ruled surface. If it is not rational, 
consider the Albanese map
$\alpha: \tilde S\to C$. Clearly $\alpha$ is 
$\Gamma$-equivariant and the action of $\Gamma$ on $C$ is not trivial.
Let $\tilde D_1\subset \tilde D$ be an $\alpha$-horizontal
component. By Adjunction $\tilde D_1$ is an elliptic curve.
So is $C$. This proves (i).

If the action on components $S_i\subset S$ 
is not transitive, 
we have an invariant divisor $S'<S$.
Since $X$ is $G\QQ$-factorial and $\rho(X)^G=1$,
we can take $S'$ so that $-(K_X+2S')$ is nef.
This contradicts Lemma \ref{lemma-hyp-sect-0}.
So, (ii) is proved.

Now we prove (iii).
Let $\Gamma$ be the image of $G_{S_i}$ in $\Aut(S_{i})$.
By Lemma \ref{lemma-fixed-points} $\rk G_{S_i}\le \rk \Gamma+1$.
If $S_i$ is rational, then we get the assertion by Theorem \ref{theorem-cr2}.
Assume that $S_i$ is birationally ruled surface over an elliptic curve.
As above, let $\tilde S_i\to S_i$ be the composition of the normalization and 
the minimal resolution, and let $\alpha: \tilde S_i\to C$ be the Albanese map. 
Then $\Gamma$ acts faithfully on $\tilde S_i$ and $\alpha$ is $\Gamma$-equivariant.
Thus we have a homomorphism $\alpha_* : \Gamma \to \Aut (C)$.
Here $\rk \Gamma \le\rk \alpha_*(\Gamma)+1+\delta_{p,2}$.
Note that $\alpha_*(\Gamma)$ is a $p$-elementary subgroup 
of the automorphism group of an elliptic curve.
Hence, $\rk \alpha_*(\Gamma)\le 2$. This implies (iii).

It remains to prove (iv).
Assume that $S$ is irreducible. By (iii)
the surface $S$ is not rational.
So, $S$ is birational to a ruled surface over an elliptic curve.
By Lemma \ref{lemma-fixed-points} the group $G$ acts 
on $S$  faithfully.
Hence, in the above notation, 
$\rk G=\rk \Gamma \le \rk \alpha_*(\Gamma)+1+\delta_{p,2}\le 3+\delta_{p,2}$, 
a contradiction.
\end{proof}

\section{Proof of Theorem \ref{theorem-cr3}}

\begin{pusto}
\label{notation-3}
In this section we prove Theorem \ref{theorem-cr3}.
As in \xref{Main-assumption} we assume that 
$X$ is a $G\QQ$-Fano threefold, where $G$ be a $p$-elementary subgroup
of $\Aut(X)$.
\end{pusto}

First we consider the case where $X$ non-Gorenstein, i.e., it has at least one 
point of index $>1$.
\begin{proposition}
\label{proposition-non-Gorenctein}
Let $G$ be a $p$-elementary group and 
let $X$ be a non-Gorenstein $G\QQ$-Fano threefold. Then 
\[
\rk G\le 
\begin{cases}
7& \text{if $p=2$,}
\\
5& \text{if $p=3$,}
\\
4& \text{if $p=5$, $7$, $11$, $13$,}
\\
3 & \text{if $p\ge 17$.}
\end{cases}
\] 
\end{proposition}
\begin{proof}
Let $P_1$ be a point of index $r>1$ and let $P_1,\dots,P_l$
be its $G$-orbit. 
Here $l=p^{t}$ for some $t$ with $t\ge s-\delta_{2,p}-3$, where $s=\rk G$ (see Lemma \ref{lemma-fixed-point}). 
By the orbifold Riemann-Roch formula \cite{Reid-YPG1987} and
a form of Bogomolov-Miyaoka 
inequality \cite{Kawamata-1992bF}, \cite{KMMT-2000} we have
\begin{equation*}
\label{eq-Kawamata-Bogomolov} 
\sum \left( r_{P_i}-\frac1{r_{P_i}} \right)<24.
\end{equation*}
Since $r_i-1/r_i\ge 3/2$, we have $3l/2<24$ and so 
\[
p^{s-\delta_{2,p}-3}\le l<16.
\]
This gives us the desired inequality.
\end{proof}

From now on we assume that our $G\QQ$-Fano threefold $X$ is 
Gorenstein, i.e., $K_X$ is a Cartier divisor.
Recall (see, e.g., \cite{Iskovskikh-Prokhorov-1999}) 
that the Picard group of a Fano variety $X$ with at worst 
(log) terminal singularities is a torsion free finitely generated 
abelian group ($\simeq H^2(X,\ZZ)$). Then we can define the \textit{Fano index} of 
$X$ as the maximal positive integer that divides $-K_X$ in $\Pic(X)$.

\begin{proposition-definition}[see, e.g., \cite{Iskovskikh-Prokhorov-1999}]
Let $X$ be a Fano threefold with at worst terminal Gorenstein
singularities. 
The positive integer $-K_X^3$ is called the \emph{degree} of $X$.
We can write $-K_X^3=2g-2$, where 
$g$ is an integer $\ge 2$ called the \emph{genus} of $X$.
Then $\dim |{-}K_X|=g+1\ge 3$.
\end{proposition-definition}

\begin{corollary-notation}
\label{corollary-notation}
In notation \xref{notation-3} the linear system
$|{-}K_X|$ is not empty, so there exists a $G$-invariant 
member $S\in |{-}K_X|$. Write $S=\sum _{i=1}^N S_i$,
where $S_i$ are irreducible components.
\end{corollary-notation}

\begin{theorem}[{\cite{Namikawa-1997}}]
\label{theorem-Namikawa}
Let $X$ be a Fano threefold with terminal Gorenstein
singularities. Then $X$ is smoothable, that is,
there is a flat family $X_t$ 
such that $X_0\simeq X$ and a general member $X_t$ is 
a smooth Fano threefold of the same 
degree, Fano index and Picard number.
Furthermore, the number of singular points is bounded as follows:
\begin{equation}
\label{eq-number-singular-points}
|\Sing(X)|\le 20-\rho(X_t)+h^{1,2}(X_t).
\end{equation}
where $h^{1,2}(X_t)$ is the Hodge number. 
\end{theorem}

Combining the above theorem with the classification of 
smooth Fano threefolds
\cite{Iskovskikh-1980-Anticanonical}, \cite{Mori1981-82}
(see also \cite{Iskovskikh-Prokhorov-1999}) we get the following

\begin{theorem}
\label{classification-Fano}
Let $X$ be a Fano threefold with at worst terminal Gorenstein
singularities and let $X_t$ be its smoothing. Let $g$ and $q$ be the genus 
and Fano index of $X$, respectively. 
\begin{enumerate}
\item
$q\le 4$.
\item
If $q=4$, then $X\simeq \PP^3$.
\item
If $q=3$, then $X$ is a quadric in $\PP^4$ \textup(with $\dim \Sing(X)\le 0$\textup).
\item
If $q=2$, then $\rho(X)\le 3$ 
and $-K_X^3=8d$, where $1\le d\le 7$. Moreover $\rho(X)=1$ if and only if
$d\le 5$.
\item
If $q=1$ and $\rho(X)=1$, then there are the following possibilities:
{\rm \begin{center}
\begin{tabular}{|c|cccccccccc|}
\hline
&&&&&&&&&&
\\[-5pt]
$g$&2&3&4&5&6&7&8&9&10&12
\\[5pt]
\hline
&&&&&&&&&&
\\[-5pt]
$h^{1,2}(X_t)$&52&30&20&14&10&7&5&3&2&0
\\[5pt]
\hline
\end{tabular}
\end{center}}
\end{enumerate}
\end{theorem}

\begin{lemma}
\label{lemma-base-point-free}
Let $G$ be a $p$-elementary group and 
let $X$ be a Gorenstein $G\QQ$-Fano threefold. If 
the linear system $|{-}K_X|$ is not base point free, then $\rk G\le 3+\delta_{p,2}$. 
\end{lemma}
\begin{proof}
Assume that $\Bs |{-}K_X|\neq \emptyset$.
Clearly, $\Bs |{-}K_X|$ is $G$-invariant.
By 
\cite{Iskovskikh-1980-Anticanonical}, \cite{Shin1989}
$\Bs |{-}K_X|$ is either a single point or a smooth rational curve.
In the first case the assertion 
immediately follows by Lemma \ref{lemma-fixed-point}.
In the second case $G$ acts on the curve $C=\Bs |{-}K_X|$.
Since $C\simeq \PP^1$, 
the assertion 
follows by Lemma \ref{lemma-fixed-points}.
\end{proof}

\begin{proposition}
\label{proposition-Gorenctein-p}
Let $G$ be a $p$-elementary group, where $p\ge 5$, and 
let $X$ be a Gorenstein $G\QQ$-Fano threefold. Then 
\[
\rk G\le 
\begin{cases}
4&\text{if $p=5,\, 7,\,11,\,13$,}
\\
3&\text{if $p\ge 17$.}
\end{cases}
\]
\end{proposition}

\begin{proof}
Assume that the above inequality does not hold.
We use the notation of \ref{corollary-notation}.
In particular, $N$ denotes the number of components of $S=\sum S_i\in |-K_X|$.
By Lemma \ref{lemma-hyp-sect} $N=p^l$, where $l\ge 1$. 
Hence $p$ divides $-K_X^3=2g-2=\bigl((-K_X)^2\cdot S_i\bigr)N$.
First we claim that $\rho(X)=1$.
Indeed, if $\rho(X)>1$, then the natural representation of $G$ on $\Pic_\QQ(X):=\Pic(X)\otimes \QQ$
is decomposed as $\Pic_\QQ(X)=V_1\oplus V$, where
$V_1$ is a trivial subrepresentation generated by the class of $-K_X$ 
and $V$ is a subrepresentation such that $V^G=0$.
Since $G$ is a $p$-elementary group, 
$\dim V\ge p-1$. 
Hence, $\rho(X)\ge p\ge 5$ and by the classification \cite{Mori1981-82} we have two possibilities:
\begin{itemize}
\item 
$-K_X^3=6(11-\rho(X))$,\quad $5\le \rho(X)\le 10$, or
\item 
$-K_X^3=28$,\quad $\rho(X)=5$.
\end{itemize}
In the last case $p=5$, so $-K_X^3\not\equiv 0\mod p$, a contradiction.
In the first case $p$ divides $-K_X^3$ only if $p=5$. 
Then $\rho(X)=6$. So, $\dim V=5$ and $V^G\neq 0$. Again we get a contradiction.

Therefore, $\rho(X)=1$. Let $q$ be the Fano index of $X$. 
We claim that $X$ is singular.
Indeed, otherwise all the $S_i$ are Cartier divisors. Then
$-K_X=NS_1$, where $N\ge p$, and so $q\ge 5$.
This contradicts (i) of Theorem \ref{classification-Fano}.
Hence $X$ is singular. By Lemma \ref{lemma-fixed-point} and our assumption we have
$|\Sing(X)|\ge p$. In particular, $q\le 2$ (see Theorem \ref{classification-Fano}).
If $q=1$, then by Theorem \ref{classification-Fano} either $2\le g\le 10$ or
$g=12$. 
Thus $N=p$ and we get the following possibilities:
$(p,g)=$ $(5,6)$, $(7,8)$, or $(11,12)$.
Moreover, $(-K_X)^2\cdot S_i= (2g-2)/N=2$.
Therefore, the restriction $|{-}K_X||_{S_i}$ of the 
(base point free) anticanonical linear system defines 
either an isomorphism to a quadric $S_i\to Q\subset \PP^3$ or
a double cover $S_i\to \PP^2$.
In both cases the image is rational, so 
we get a map $G_i \to \Cr_2(\CC)$ whose kernel is of rank $\le 1$ 
by Lemma \ref{lemma-fixed-points} and because $p>2$.
Then by Theorem 
\ref{theorem-cr2} \ 
$\rk G_{S_i}\le 3$.
Hence, $\rk G\le 4$ which contradicts our assumption.

Finally, consider the case $q=2$. Then $-K_X=2H$
for some ample Cartier divisor $H$ and $d:=H^3\le 7$.
Therefore, $N S_i\cdot H^2=S\cdot H^2=2d$.
Since $\rho(X)=1$, by Theorem \ref{classification-Fano}
we get $p=d=5$. Then we apply \eqref{eq-number-singular-points}.
In this case, $h^{1,2}(X_t)=0$  (see \cite{Iskovskikh-Prokhorov-1999}). So,
$|\Sing(X)|\le 19$. On the other hand, 
$|\Sing(X)|\ge 25$ by Lemma \ref{lemma-fixed-point} and our assumption.
The contradiction proves the proposition.
\end{proof}

We need the following result which is a very weak form of much more general Shokurov's 
toric conjecture \cite{McKernan-2001}, \cite{Prokhorov-2003e}.
\begin{lemma}
\label{lemma-toric-conjecture}
Let $V$ be a smooth Fano threefold and let 
$D\in |{-}K_V|$ be a divisor such that the pair
$(V,D)$ is LC. 
Then $D$ has at most $3+\rho(V)$ irreducible components.
\end{lemma}

\begin{proof}
Write $D=\sum_{i=1}^n D_i$.
If $\rho(V)=1$, then all the $D_i$ are linearly proportional:
$D_i\sim n_i H$, where $H$ is an ample generator of $\Pic(V)$.
Then $-K_V\sim \sum n_i H$ and 
by Theorem \ref{classification-Fano} we have $\sum n_i=q\le 4$.

If $V$ is a blowup of a curve on another smooth Fano threefold $W$,
then we can proceed by induction replacing $V$ with $W$.
Thus we assume that $V$ cannot be obtained by blowing up 
of a curve on another smooth Fano threefold. In this situation $V$  is  
called \textit{primitive} (\cite{Mori1983}). 
According to \cite[Th. 1.6]{Mori1983} we have $\rho(V)\le 3$
and $V$ has a conic bundle structure $f: V\to Z$, where $Z\simeq \PP^2$
(resp. $Z\simeq \PP^1\times \PP^1$)
if $\rho(V)=2$ (resp. $\rho(V)=3$). Let $\ell$ be a general fiber.
Then $2=-K_V\cdot \ell=\sum D_i\cdot \ell$.
Hence $D$ has at most two $f$-horizontal components 
and at least $n-2$ vertical ones. 
Now let $h: V\to W$ be an extremal contraction other than $f$
and let $\ell'$
be any curve in a non-trivial fiber of $h$.
For any $f$-vertical component
$D_i\subset D$ we have $D_i=f^{-1}(\Gamma_i)$, where  $\Gamma_i\subset Z$ is a curve, 
so $D_i\cdot \ell'=\Gamma_i \cdot f_*\ell'\ge 0$.
If $\rho(V)=2$, then $D_i\cdot l'\ge 1$. Hence, $-K_V\cdot \ell'\ge n-2$.
On the other hand, $-K_V\cdot \ell'\le 3$ (see \cite[\S 3]{Mori1983}).
This immediately gives us $n\le 5$ as claimed.
Finally consider the case $\rho(V)=3$.
Assume that $n\ge 7$. Then we can take $h$ so that  
$\ell'$ meets at least three $f$-vertical components,
say $D_1$, $D_2$, $D_3$. As above, $-K_V\cdot \ell'\ge 3$
and by the classification of extremal rays (see \cite[\S 3]{Mori1983})
$h$ is a del Pezzo fibration. 
This contradicts our assumption $\rho(V)=3$.
\end{proof}

\begin{proposition}
\label{proposition-Gorenctein-2}
 Let $G$ be a $2$-elementary group and 
let $X$ be a Gorenstein $G\QQ$-Fano threefold. Then $\rk G\le 7$.
\end{proposition}
\begin{proof}
Assume that $\rk G\ge 8$. 
By Lemma \ref{lemma-hyp-sect} we have 
$\rk G_{S_i}\le 5$. Hence, $N\ge 8$.
If $X$ is smooth, then 
by Lemma \ref{lemma-toric-conjecture} we have $\rho(X)\ge 5$.
If furthermore $X\simeq Y\times \PP^1$, where $Y$ is a del Pezzo surface, then
the projection $X\to Y$ must be $G$-equivariant.
This contradicts $\rho(X)^G=1$. 
Therefore, $\rho(X)=5$ and $-K_X^3=28$ or $36$ (see \cite{Mori1981-82}).
On the other hand, $-K_X^3$ is divisible by $N$, a contradiction.


Thus $X$ is singular.
Assume that $|\Sing(X)|\ge 32$.
Then for a smoothing $X_t$ of $X$
by \eqref{eq-number-singular-points} we have 
$h^{1,2}(X_t)\ge 13$. Since $N$ divides $-K_X^3=-K_{X_t}^3$,
using the classification of Fano threefolds 
\cite{Iskovskikh-1980-Anticanonical}, \cite{Mori1981-82}
(see also \cite{Iskovskikh-Prokhorov-1999}) we get:
\[
\rho(X)=1,\quad -K_X^3=8,\quad  N=8,\quad  |\Sing(X)|=32.
\]
Consider the representation of $G$ on $H^0(X,-K_X)$.
Since 
\[
 7=\dim H^0(X,-K_X)< \rk G, 
\]
this representation is not faithful (otherwise $G$ is contained in 
a maximal torus of $\GL(H^0(X,-K_X))=\GL_7(\CC)$).
Therefore, the linear system $|{-}K_X|$ is not very ample.
On the other hand, $|{-}K_X|$ is base point free (see 
Lemma \ref{lemma-base-point-free}). Hence $|{-}K_X|$ defines a double cover 
$X\to Y\subset \PP^{6}$ \cite{Iskovskikh-1980-Anticanonical}.
Here $Y$ is a variety of degree $4$ in $\PP^6$,
a variety of minimal degree. 
If $Y$ is smooth, then according to the Enriques theorem 
(see, e.g., \cite[Th. 3.11]{Iskovskikh-1980-Anticanonical})
$Y$ is a rational scroll
$\PP_{\PP^1}(\EEE)$, where 
$\EEE$ is a rank $3$ vector bundle on $\PP^1$.
Then $X$ has a $G$-equivariant projection to a curve.
This contradicts $\rho(X)^G=1$.
Hence $Y$ is singular. In this case, $Y$ is a cone (again by the the Enriques theorem 
\cite[Th. 3.11]{Iskovskikh-1980-Anticanonical}).
If its vertex $O\in Y$ is zero-dimensional, then  $\dim T_{O,Y}=6$.
On the other hand, $X$ has only hypersurface singularities (see \ref{Terminal-singularities}).
Therefore the double cover $X\to Y$ is not \'etale over $O$ and so
$G$ has a fixed point 
on $X$. This contradicts Lemma \ref{lemma-fixed-point}.
Thus $Y$ is a cone over a rational normal curve of degree $4$
with vertex along a line. Then $X$ cannot have isolated singularities, a contradiction.

Therefore, $|\Sing(X)|< 32$.
Then for any point $P\in\Sing(X)$ by Lemma \ref{lemma-fixed-point} 
we have $\rk G_P \ge 4$. Hence the orbit of $P$ contains $16$ elements
and coincides with $\Sing(X)$, i.e. the action of $G$ on $\Sing(X)$ is transitive.
Since $S\cap \Sing(X)\neq \emptyset$, we have $\Sing(X)\subset S$.
On the other hand, our choice  of $S$ in \ref {Main-assumption} is not unique:
there is a basis $s^{(1)}$, \dots, $s^{(g+2)}\in H^0(X,-K_X)$
consisting of eigensections. This basis gives us $G$-invariant divisors 
$S^{(1)}$, \dots, $S^{(g+2)}$ generating $|-K_X|$.
By the above $\Sing(X)\subset S^{(i)}$ for all $i$.
Thus $\Sing(X)\subset \cap S^{(i)}=\Bs |-K_X|$.
This contradicts  Lemma \ref{lemma-base-point-free}.
Proposition \ref{proposition-Gorenctein-2} is proved.
\end{proof}

\begin{proposition}
\label{proposition-Gorenctein-3}
Let $G$ be an $3$-elementary group and 
let $X$ be a Gorenstein $G\QQ$-Fano threefold. Then $\rk G\le 5$.
\end{proposition}
\begin{proof}
Assume that $\rk G\ge 6$. 
By Lemma \ref{lemma-hyp-sect} we have 
$\rk G_{S_i}\le 5$.
Hence, $N\ge 9$.
If $X$ is smooth, then by Lemma \ref{lemma-toric-conjecture} we have $\rho(X)\ge 6$ and so
$X\simeq Y\times \PP^1$, where $Y$ is a del Pezzo surface \cite{Mori1981-82}.
Then the projection $X\to Y$ must be $G$-equivariant.
This contradicts $\rho(X)^G=1$.
Therefore, $X$ is singular.
By Lemma \ref{lemma-fixed-point} $|\Sing(X)|\ge 3^{6-3}=27$. 
Hence, for a smoothing $X_t$ of $X$
by \eqref{eq-number-singular-points} we have $h^{1,2}(X_t)\ge 7+\rho(X)$. 
Recall that $N$ divides $-K_X^3=-K_{X_t}^3$.
Then we use the classification of smooth Fano threefolds 
\cite{Iskovskikh-1980-Anticanonical}, \cite{Mori1981-82}
and get a contradiction.
\end{proof}

Now Theorem \ref{theorem-cr3} is a consequence of Propositions 
\ref{proposition-non-Gorenctein},
\ref{proposition-Gorenctein-p},
\ref{proposition-Gorenctein-2}, and
\ref{proposition-Gorenctein-3}.

\def\cprime{$'$}
  \def\mathbb#1{\mathbf#1}\def\polhk#1{\setbox0=\hbox{#1}{\ooalign{\hidewidth
  \lower1.5ex\hbox{`}\hidewidth\crcr\unhbox0}}}

\end{document}